\documentclass[letterpaper, 10 pt, conference]{ieeeconf}

\IEEEoverridecommandlockouts 
\usepackage{graphicx}
\usepackage{microtype}
\usepackage{lmodern}
\usepackage[most]{tcolorbox}

\newcommand{\behcetFull}{Beh\c{c}et~A\c{c}\i kme\c{s}e}

\newcommand{\bma}[1]{\left[\begin{array}{#1}}
\newcommand{\ema}{\end{array}\right]}

\newcommand{\minimize}[1]{\underset{#1}{\operatorname{minimize}}}

\newcommand{\dtlos}{{\small\textsc{DT-LoS}}}
\newcommand{\ctlos}{{\small\textsc{CT-LoS}}}
\newcommand{\ctscvx}{{\small\textsc{CT-SCvx}}}
\newcommand{\st}{\mathrm{subject}\ \mathrm{to}}
\definecolor{steelblue}{HTML}{8CC9F7}
\definecolor{green}{HTML}{8AEEA8}

\newtcolorbox{mybox}{
  enhanced,
  colback=steelblue!50!white,
  colframe=steelblue,
  fonttitle=\bfseries,
  boxrule=0pt,  % Adjust this value for border width
  drop shadow
}

\newtcolorbox{myboxyellow}{
  enhanced,
  colback=green!50!white,
  colframe=yellow,
  fonttitle=\bfseries,
  boxrule=0pt,  % Adjust this value for border width
  drop shadow,
  boxsep=0pt,  % space between frame and content
  left=2mm,  % left margin
  right=2mm,  % right margin
  top=2mm,  % top margin
  bottom=2mm  % bottom margin
}

\usepackage{algorithm}
\usepackage{algpseudocode}
\algrenewcommand{\algorithmiccomment}[1]{\hfill\textcolor{black!65!white}{$\triangleright$ #1}}

\usepackage{setspace}
\usepackage{varwidth} 
\algnewcommand{\LineComment}[1]{\State \textcolor{black!65!white}{$\triangleright$ #1}}
\newcommand\copyrighttext{
  \footnotesize \textcopyright\ 2025 IEEE.  Personal use of this material is permitted.  Permission from IEEE must be obtained for all other uses, in any current or future media, including reprinting/republishing this material for advertising or promotional purposes, creating new collective works, for resale or redistribution to servers or lists, or reuse of any copyrighted component of this work in other works.}
\newcommand\copyrightnotice{
\begin{tikzpicture}[remember picture,overlay]
\node[anchor=south,yshift=10pt] at (current page.south) {\fbox{\parbox{\dimexpr\textwidth-\fboxsep-\fboxrule\relax}{\copyrighttext}}};
\end{tikzpicture}}

\usepackage{amsmath, amssymb, amsfonts} 
\allowdisplaybreaks
\usepackage{tabularray}
\usepackage{xcolor}
\usepackage{subcaption}
\usepackage{caption}
\usepackage{siunitx}
\usepackage{pifont}% http://ctan.org/pkg/pifont
\usepackage{hyperref}
\newcommand{\cmark}{\textcolor{green}{\ding{51}}}%
\newcommand{\xmark}{\textcolor{red}{\ding{55}}}%
\usepackage{bm}
\usepackage{float}

\title{\LARGE \bf Continuous-Time Line-of-Sight Constrained Trajectory Planning for 6-Degree of Freedom Systems}

\author{Christopher R. Hayner$^{1}$, John M. Carson III$^2$, \behcetFull$^1$, and Karen Leung$^1$ % <-this % stops a space
\thanks{$^1$ Dept. of Aeronautics and Astronautics, University of Washington, Seattle, WA, USA. Email: {\tt\small haynec@uw.edu}}
\thanks{$^2$ NASA Johnson Space Center, Houston, TX, USA}
\thanks{This work was supported by a NASA Space Technology Graduate Research Opportunity and the Office of Naval Research under grant N00014-17-1-2433.}
\thanks{The authors would like to acknowledge Natalia Pavlasek, Griffin Norris, Samuel Buckner, and Purnanand Elango for their many helpful discussions and support throughout this work.}
\thanks{The code and implementation details for this work can be found \href{https://haynec.github.io/papers/los}{https://haynec.github.io/papers/los}.}}

\begin{document}

\maketitle
\copyrightnotice
\thispagestyle{empty}
\pagestyle{empty}

\begin{abstract}
Perception algorithms are ubiquitous in modern autonomy stacks, providing necessary environmental information to operate in the real world.
Many of these algorithms depend on the visibility of keypoints, which must remain within the robot’s line-of-sight (LoS), for reliable operation. 
This paper tackles the challenge of maintaining LoS on such keypoints during robot movement. 
We propose a novel method that addresses these issues by ensuring applicability to various sensor footprints, adaptability to arbitrary nonlinear system dynamics, and constant enforcement of LoS throughout the robot's path. 
Our experiments show that the proposed approach achieves significantly reduced LoS violation and runtime compared to existing state-of-the-art methods in several representative and challenging scenarios.

\end{abstract}

\section{Introduction}
\begin{figure}[t]
    \centering
    \includegraphics[width=\linewidth]{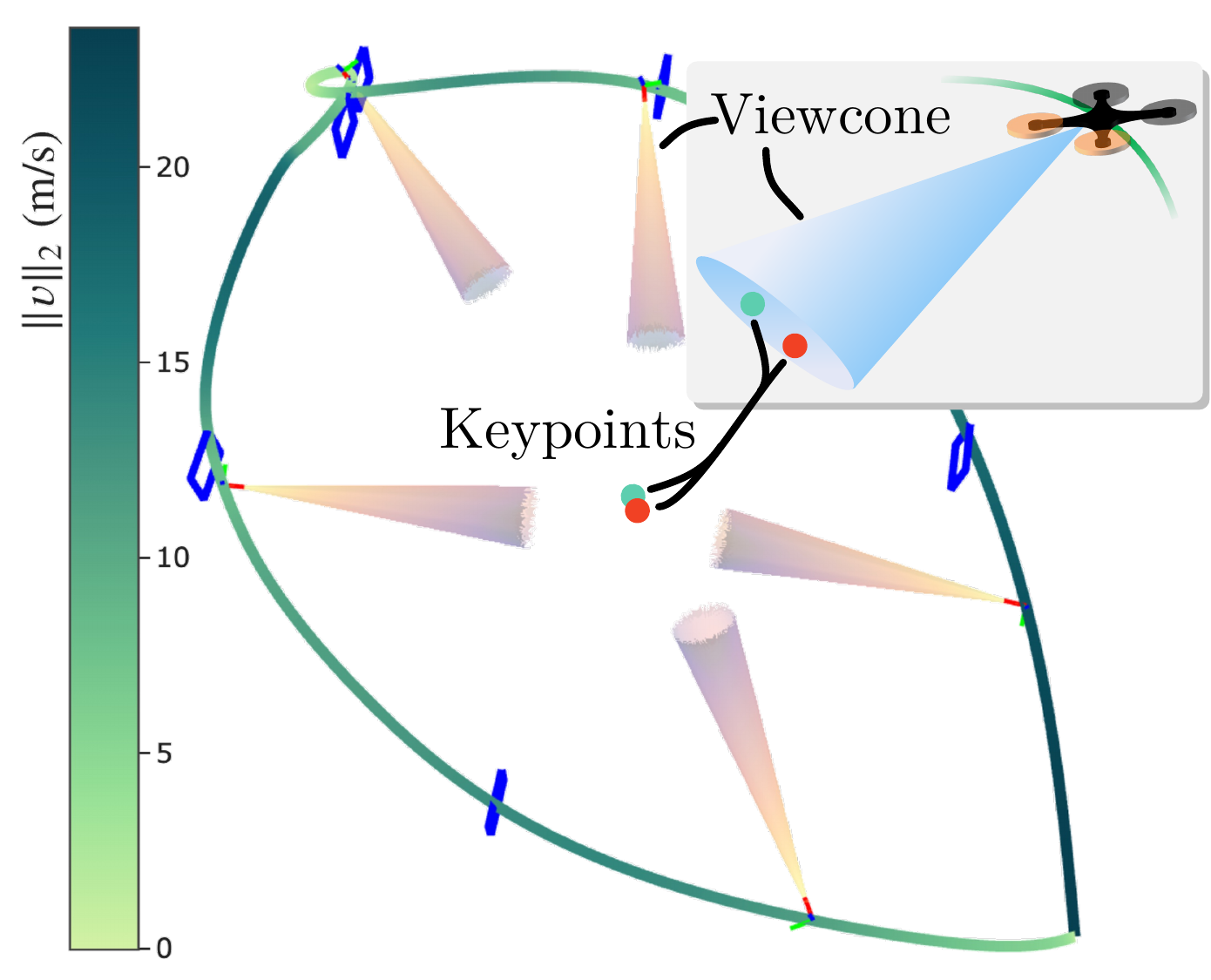}
    \caption{LoS Guidance: A simulation result of our LoS guidance method showing the pose of the quadrotor and view cone at several points along the trajectory. Throughout the trajectory, the keypoints are within LoS.}
    \label{fig:main}
    \vspace{-12pt}
\end{figure}

Perception algorithms are essential for enabling robots to perform relative navigation and tracking tasks, allowing them to gather information about internal and external states for \textit{safety-critical} tasks (\textit{e.g.}, pose estimation, obstacle avoidance), without relying on external sources (\textit{e.g.}, GNSS, motion capture, cooperative agents).
Given the increasing interest in deploying robotics in real-world unknown environments, there is a heavy reliance on perception algorithms for the aforementioned safety-critical tasks, perception algorithms must be reliable.
A central requirement for many of these algorithms, which depend on observing keypoints or points of interest, is that these keypoints must always remain visible to the robot.
For example, landmarks within visual-inertial odometry, \cite{Huang2019-so, Scaramuzza2011-bg}, a subject being filmed in drone-enabled cinematography, \cite{Huang2018-pi, Alcantara2021cine}, or prospective landing sites and hazards for a planetary lander, \cite{Hayner2023halo, Carson2019-va, Johnson2008trn}. 
If keypoints do not remain visible while the robot is moving, perception performance may suffer likely resulting in the task being considered a failure. 
In this paper, we consider the problem of trajectory optimization under LoS constraints and refer to it as \textit{LoS guidance}, depicted in Fig.~\ref{fig:main}.

State-of-the-art methods for LoS guidance in these applications do not satisfy all the key requirements of an effective, realizable, and generalizable LoS guidance method:
\begin{itemize}
    \item[] \textbf{R1.} applicability to different sensor footprints, 
    \item[] \textbf{R2.} adaptability to arbitrary nonlinear dynamics, and
    \item[] \textbf{R3.} maintaining LoS on keypoints throughout the entire trajectory.
\end{itemize}
We present a method of LoS guidance that addresses these requirements.

\noindent\textbf{Contributions.} The contributions of this paper are three-fold: \textbf{1.} a sensor-footprint-agnostic LoS constraint, \textbf{2.} a computationally tractable formulation of a continuous-time LoS guidance methodology, \ctlos , as well as a baseline discrete-time LoS \dtlos\ guidance methodology, and finally \textbf{3.} a demonstration of the efficacy of \ctlos\ and a comparison against \dtlos\ in several keypoint tracking problems inspired by challenges encountered in relative navigation and drone-enabled cinematography.

\noindent\textbf{Organization.} This paper is organized into the following sections. Section~\ref{sec:lit_rev} discusses related work in the field.
Section~\ref{sec:gen_problem} introduces the general problem as well as the notation used throughout this paper. Section~\ref{sec:6dof} introduces the rigid-body dynamics and six-degree-of-freedom LoS constraint. \ctlos\ is introduced as the solution method in Section~\ref{sec:scp}. 
The experiments and baseline comparison are introduced in Section~\ref{sec:experiments}.
Results from these experiments are discussed in Section~\ref{sec:discuss}. Finally, conclusions are discussed in Section~\ref{sec:conclusion}.
\noindent\textbf{Terminology:} We use the term \textit{six-degree-of-freedom} when referring to the three translational modes and three rotational modes of a rigid body.

\section{Literature Review} \label{sec:lit_rev}
In this section, we offer a brief taxonomy of nonlinear trajectory optimization techniques with LoS constraints, seen in Table~\ref{tab:traj_class}, and justification for our formulation and methodology. 
\begin{table}[h]
    \small
    \centering
    \begin{minipage}{0.49\textwidth}
        \centering
        \begin{tblr}{
            colspec={|X[c]|X[c,1.5cm]|X[c,1.5cm]|X[c,1.5cm]|},
            width=\linewidth
        }
             \hline
             Method & \textbf{R1} & \textbf{R2} & \textbf{R3} \\
             \hline
             \cite{Zhou2021-nl, Tord2022panther, Spasojevic2020-vc, Murali2019-qf} & \xmark & \xmark & \xmark \\
             \hline
             \cite{Falanga2018vio, reynolds2020dqg} & \xmark & \cmark & \xmark \\
             \hline
        \end{tblr}
    \end{minipage}
    \begin{minipage}{0.49\textwidth}
        \centering
        \begin{tblr}{
            colspec={|X[c]|X[c,1.5cm]|X[c,1.5cm]|X[c,1.5cm]|},
            width=\linewidth
        }
             \hline
             Baseline & \cmark & \cmark & \xmark \\
             \hline
             Proposed & \cmark & \cmark & \cmark\\
             \hline
        \end{tblr}
    \end{minipage}
    \caption{Classification of Related Work}
    \label{tab:traj_class}
    \vspace{-10pt}
\end{table}
\normalsize

Pseudospectral trajectory planning, as discussed in Section 9.7 \cite{Kelly2017-tt}, is a method that parameterizes state and control profiles using polynomial splines, with coefficients or knot points as decision variables in the Nonlinear Program (NLP). 
For quadrotors, these methods often use the differentially flat model from \cite{Mellinger2011-qw}, optimizing position and yaw angle trajectories either independently or jointly. 
They are effective at generating online trajectories \cite{Richter2016-hk, Ding2019-cm}.

For LoS guidance, \cite{Zhou2021-nl} proposed a pseudospectral method following a two-step decoupled approach, maximizing information gain by adjusting yaw angles to include the regions in LoS with the highest uncertainty. 
However, this method does not strictly constrain points of interest to be within the sensor’s LoS. 
Similarly, \cite{Spasojevic2020-vc} and \cite{Murali2019-qf} use a decoupled approach but enforce a symmetric second-order cone constraint to keep keypoints within the LoS. 
However, the decoupled approach results in greedy solutions where positions are optimized without consideration of LoS constraints. This leads to solutions that unnecessarily sacrifice optimality with respect to their objective.

To address this, \cite{Tord2022panther} jointly optimizes position and yaw angle.
This allows for overall less conservative solutions when compared to the decoupled approach.
However, vehicle dynamic feasibility is limited to convex min-max constraints on the flat states and their derivatives, and cannot be generalized to arbitrary nonlinear dynamics.

In contrast to pseudospectral methods, \cite{Falanga2018vio} adopts a multiple-shooting method (see Section 9.6 in \cite{Kelly2017-tt}) and directly imposes the nonlinear quadrotor dynamics constraint with no assumption of differential flatness at each discrete node. 
However, keypoint inclusion in LoS, modeled by a non-symmetric $L_{\infty}$-norm cone, restricting the choice of sensor footprint, is only promoted in the cost function rather than directly enforced as a constraint. 
This results in insufficient satisfaction of LoS constraints.

All of the aforementioned methods leverage software such as ACADO \cite{Houska2011-ju} to solve the NLP.
These packages call underlying NLP solvers such as IPOPT \cite{Wachter2006-uo} or convex approximation algorithms such as sequential quadratic programming \cite{Boggs1995-av}. A limitation of relying on these solvers is they lack convergence guarantees, rely on second-order information, and they cannot enforce constraints in continuous time.
In contrast, successive convexification-based methods (SCvx) has gained significant popularity in trajectory optimization, \cite{Malyuta2022csm, Mao2017-zj}, as they have convergence guarantees \cite{Bonalli2019-jj}, only require first-order information\footnote{when using a first-order solver, such as \cite{Yu2022-mb}}, and, recently, have guarantees on continuous-time constraint satisfaction \cite{Elango2024ctcs}.
Sequential Convex Programming (SCP) has been used to solve the LoS guidance problem in a planetary landing context \cite{reynolds2020dqg}, a symmetric $L_2$ norm cone constraint is enforced at each discrete node to ensure a single landing site is within LoS.
However, these methods only consider a single norm type and static keypoint, meaning they do not consider different sensor footprints or dynamic keypoints.
Finally, a major limitation of all the LoS guidance methods discussed in this section is that they only enforce the LoS constraint at discrete nodes along the trajectory, rather than in continuous-time. Therefore, they lack constraint satisfaction in between nodes, rendering them hazardous in safety-critical situations.

\section{General Problem Formulation}\label{sec:gen_problem}
In this section, we introduce the notation used throughout the paper in Section~\ref{sec:not} and the general continuous-time LoS guidance problem formulation in Section~\ref{sec:gen_form}.

\subsection{Notation}
\label{sec:not}
We use the following notation in the remainder of this paper. 
A quantity expressed in the frame $\mathcal{A}$ is denoted by the subscript $\Box_{\mathcal{A}}$. 

To parameterize the attitude of frame $\mathcal{B}$ with respect to frame $\mathcal{A}$, the unit quaternion, $q_{\mathcal{A} \to \mathcal{B}} \in \mathcal{S}^3$ where $\mathcal{S}^3\subset\mathbb{R}^4$ is the unit 3-sphere, is used. 
The scalar-first convention is used for unit quaternion as follows $q = \begin{bmatrix} q_w & q_x & q_y & q_z \end{bmatrix}^\top$. 

\subsection{Continuous-Time Formulation}
\label{sec:gen_form}
We consider a keypoint to be in the LoS of an exteroceptive sensor if it resides within a view cone, as seen in Fig~\ref{fig:cone}. 
We frame the LoS guidance problem as a general, continuous, free-final-time trajectory planning problem~\eqref{eq:gen_prob}.

\begin{mybox}{\textbf{General LoS Problem}}
\small
\begin{subequations}
\label{eq:gen_prob}
\begin{alignat}{2}
    &\!\minimize{x,u, t_f} &\quad& L_{f}(t_f, x(t_f),u(t_f)), \label{eq:1}\\
    &\st & &\dot{x}(t) = f(t, x(t),u(t)), \label{eq:non_con_prob_dyn}\\
    & & &p(t) \in \mathcal{K}(b(x(t)), c(x(t)), A^{\mathrm{C}}(t)), \label{eq:non_con_los}\\
    & & &g(t,x(t),u(t)) \leq 0_{n_g},\\
    & & &h(t,x(t),u(t)) = 0_{n_h},\\
    & & &x(t_i) = x_i, x(t_f) = x_f,
\end{alignat}
\end{subequations}
where appropriate, $\forall t \in [t_i,t_f]$.
\end{mybox}
\normalsize
\noindent In \eqref{eq:gen_prob} $L_
f:\mathbb{R}_+ \times \mathbb{R}^{n_x} \times \mathbb{R}^{n_u} \mapsto \mathbb{R}$ denotes a general state, control and time-dependent cost, $f: \mathbb{R}_+\times \mathbb{R}^{n_x} \times \mathbb{R}^{n_u} \mapsto \mathbb{R}^n$ denotes the system dynamics, and the set $$\mathcal{K}(b, c, A^{\mathrm{C}}) = \{a\in\mathbb{R}^n \mid \lVert A^{\mathrm{C}}(a-b)\rVert \leq c^\top (a - b)\},$$ where $A^{\mathrm{C}}: \mathbb{R}^{n\times n} \mapsto \mathbb{R}^{n\times n}, b\in\mathbb{R}^{n}, c\in\mathbb{R}^{n}$, denotes the view cone. Here a point to be kept within LoS is denoted by $a$, the position of the sensor is denoted by $b(t)$, the boresight vector is denoted by $c(t)$, and the angular frame of view is denoted by $A^{\mathrm{C}}(t)$ (for the remainder of this paper, the boresight vector and angular frame of view are kept constant).
The path constraint functions $g : \mathbb{R}_+ \times \mathbb{R}^{n_x} \times \mathbb{R}^{n_u} \mapsto \mathbb{R}^{n_g}$ denotes general inequality constraints, $h : \mathbb{R}_+ \times \mathbb{R}^{n_x} \times \mathbb{R}^{n_u} \mapsto \mathbb{R}^{n_h}$ denotes general equality constraints, and finally $x_{i} \text{ and } x_{f}$ denote initial and terminal-state constraints. We make no assumptions on the convexity of Prob.~\ref{eq:gen_prob}, nor are there assumptions of differential flatness on the dynamics in \eqref{eq:non_con_prob_dyn}.

\section{Six-Degree-of-Freedom Formulation}
\label{sec:6dof}
This section gives an overview of the six-degree-of-freedom-dynamics in Section~\ref{sec:dyn} and introduces the proposed LoS constraint in Section~\ref{sec:los}.
\subsection{Dynamics}\label{sec:dyn}
We adopt the six-degree-of-freedom (6-DoF) rigid-body dynamics similar to those used in \cite{Szmuk2019-bt, Jacquet2022vis}.
The vehicle state is defined as, $x = \begin{bmatrix} r_\mathcal{I}^\top & v_\mathcal{I}^\top & q_{\mathcal{B\to I}}^\top & \omega_\mathcal{B}^\top \end{bmatrix}^\top$, where the position resolved in the inertial frame is $r_\mathcal{I}\in\mathbb{R}^3$, the velocity resolved in the inertial frame is $v_\mathcal{I}\in\mathbb{R}^3$, the attitude of the body frame relative to the inertial frame is $q_{\mathcal{B\to I}}$, and angular rate of the body frame relative to an inertial frame is $\omega_\mathcal{B} \in \mathbb{R}^3$. 
The vehicle control vector, $u(t) = \begin{bmatrix} f_\mathcal{B}^\top & M_\mathcal{B}^\top \end{bmatrix}^\top$, consists of the thrust vector represented in the body frame, $f_\mathcal{B}\in \mathbb{R}^3$, and moment vector, $M_\mathcal{B}\in\mathbb{R}^3$. The vehicle state evolves according to the following dynamics
\small 
\begin{align*}
    % \label{eq:6dof_def}
    & \textbf{Position.} && \dot{r}_\mathcal{I}(t) = v_\mathcal{I}(t),\\
    & \textbf{Velocity.} && \dot{v}_\mathcal{I}(t) = \frac{1}{m}\left(C(q_{\mathcal{B \to I}}(t)) f_{ \mathcal{B}}(t)\right) + g_{\mathcal{I}},\\
    & \textbf{Attitude.}  && \dot{q}_{\mathcal{I}\to \mathcal{B}} = \frac{1}{2} \Omega(\omega_\mathcal{B}(t))  q_{\mathcal{I \to B}}(t),\\
    & \textbf{Angular Rate.}  && \dot{\omega}_\mathcal{B}(t) =  J_{\mathcal{B}}^{-1} \left(M_{\mathcal{B}}(t) - \left[\omega_\mathcal{B}(t)\times\right]J_{\mathcal{B}} \omega_\mathcal{B}(t) \right),
\end{align*} 
\normalsize
where $g_\mathcal{I}\in\mathbb{R}^3$ is the gravity of Earth expressed in the inertial reference frame, and $J_\mathcal{B}$ is the inertial tensor of the vehicle expressed in the body frame. 
The operator $C:\mathcal{S}^3\mapsto SO(3)$ represents the direction cosine matrix (DCM), where $SO(3)$ denotes the special orthogonal group.
For a vector $\xi \in \mathbb{R}^3$, the skew-symmetric operators $\Omega(\xi)$ and $[\xi \times]$ are defined in Section II.A. in \cite{Szmuk2018free}. Using the aforementioned definitions, the 6-DoF dynamics are

\footnotesize
\begin{multline*}
    f_{\mathrm{6DOF}}(t,x(t),u(t)) = \begin{bmatrix} \dot{r}_\mathcal{I}(t)^\top & \dot{v}_{\mathcal{I}}(t)^\top & \dot{q}_{\mathcal{I}\to\mathcal{B}}(t)^\top & \dot{\omega}_{\mathcal{B}}(t)^\top \end{bmatrix}^\top.
\end{multline*}
\normalsize

\subsection{Nonconvex Line-of-Sight Formulation}\label{sec:los}
The LoS constraint ensures that a set of keypoints remain within a view cone, $\mathcal{K}$, which models an exteroceptive sensor's frame of view. This sensor is assumed to be rigidly attached to the body of the vehicle.

As such, both the position and attitude of the sensor are solely determined by the position and attitude of the vehicle. The LoS constraint can be broken up into two parts: \textbf{1.} a nonconvex transformation component that expresses the transform from the inertial frame to the sensor frame and \textbf{2.} a convex norm cone component expressed in the sensor frame.
\subsubsection{Transformation Component}
Given the location of a keypoint represented in the inertial frame, $p_\mathcal{I}(t)$ and the vehicle state, $x(t)$, the sensor frame representation of the keypoint is
\begin{align}
    \label{eq:trans}
    p_\mathcal{S}(t) = C(q_{\mathcal{\mathcal{S}\to\mathcal{B}}})C(q_{\mathcal{\mathcal{B}\to\mathcal{I}}}(t))(p_{\mathcal{I}}(t) - r_{\mathcal{I}}(t)).
\end{align}
\subsubsection{Norm Cone Component}
\begin{figure}
    \centering
    \includegraphics[width=0.4\textwidth]{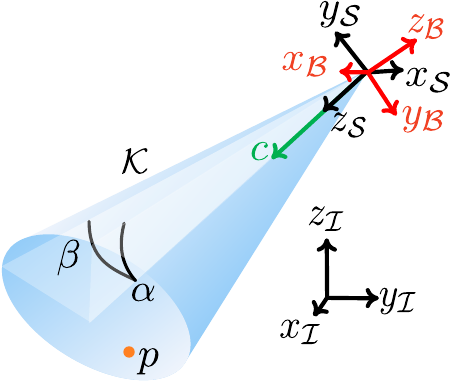}
    \caption{Norm Cone: The inertial, body, and sensor frames are denoted by the axes with subscripts $\mathcal{I,\ B}$, and $\mathcal{S}$ respectively. The keypoint is shown by the orange point $p$, the boresight vector, $c$, of the cone, $\mathcal{K}$, is represented by the green vector, and the angles of the view cone $\alpha$ and $\beta$ are visualized on the cone.}
    \label{fig:cone}
    \vspace{-10pt}
\end{figure}

Given a keypoint resolved in the sensor frame, $p_\mathcal{S}(t)$, the sensor view cone is modeled using a norm cone,
\begin{align}
    \label{eq:view_cone}
    \lVert A^{\mathrm{C}} p_\mathcal{S}(t)\rVert_\rho \leq c^\top p_\mathcal{S}(t).
\end{align}
The $z$-axis of the sensor frame is chosen to be aligned with the boresight vector, $c$, of the cone, meaning $c = \begin{bmatrix} 0 & 0 & 1 \end{bmatrix}^\top$. It follows that $A^{\mathrm{C}}$ is defined as 
\begin{align*}
    A^{\mathrm{C}} \triangleq \begin{bmatrix} \frac{1}{\tan(\alpha)} & 0 & 0 \\ 0 & \frac{1}{\tan(\beta)} & 0 \\ 0 &  0 & 0 \end{bmatrix},
\end{align*}
where $\alpha$ and $\beta$ are the angles of the field of view that are aligned with the $x_\mathcal{S}$ and $y_\mathcal{S}$ axis respectively. The norm cone is visualized in Fig~\ref{fig:cone}.
As most exteroceptive sensors have either rectangular or circular frames of view, either $\rho = \infty$ or $\rho = 2$ respectively are the appropriate choice of $\rho$-norm for \eqref{eq:view_cone}. However, this method extends to any norm. Notably, \eqref{eq:view_cone} is a convex constraint in $p_\mathcal{S}(t)$ (see Section 2.2.3 in \cite{Boyd2004cvx}).

\subsubsection{Line-of-sight Constraint}
Substituting \eqref{eq:trans} into \eqref{eq:view_cone} yields the full nonconvex LoS constraint
\begin{multline}
    \label{eq:los_constr}
    g_{\mathrm{LoS}} \triangleq \lVert A^{C} C(q_{\mathcal{S}\to\mathcal{B}})C(q_{\mathcal{B}\to\mathcal{I}}(t))(p_{\mathcal{I}}(t) - r_{\mathcal{I}}(t))\rVert_\rho - \cr c^\top C(q_{\mathcal{S}\to\mathcal{B}})C(q_{\mathcal{B}\to\mathcal{I}}(t))(p_{\mathcal{I}}(t) - r_{\mathcal{I}}(t)) \leq 0,
\end{multline}
and define it as $g_\mathrm{LoS}$ for notational brevity.

When \eqref{eq:los_constr} is satisfied, $p \in \mathcal{K}$.

\section{\ctlos\ Algorithm}\label{sec:scp}
Traditional SCvx-based methods solve the continuous-time LoS guidance problem, Prob.~\ref{eq:gen_prob}, by approximating it in discrete time and enforcing constraints at each discrete node. 
Choosing a sparse discretization grid can lead to large inter-nodal constraint violations while increasing nodes increases the computational costs.

The \ctscvx\ algorithm by \cite{Elango2024ctcs} solves the problem directly by integrating constraints over the entire trajectory, avoiding this approximation.

\begin{figure}[h!]
    \centering
    \includegraphics[width=0.45\textwidth]{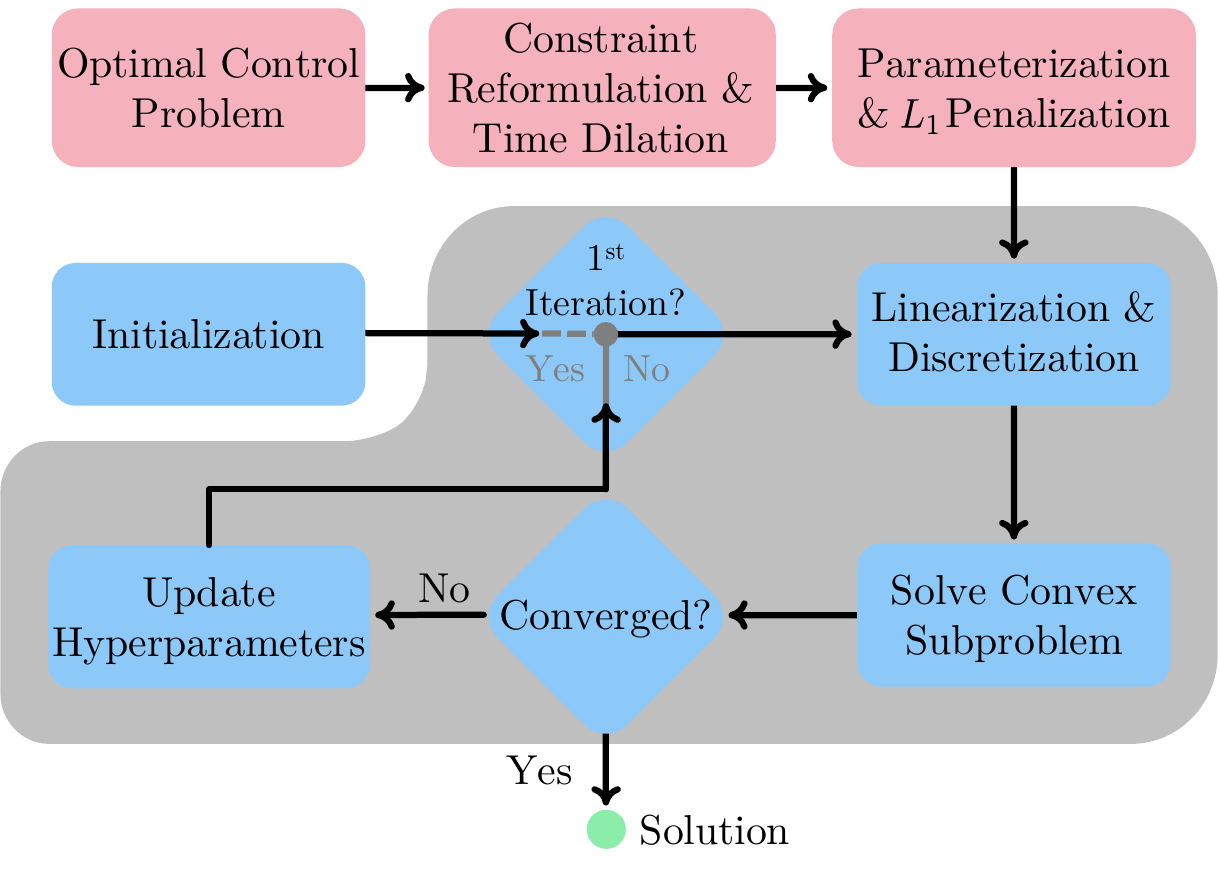}
    \caption{The \ctscvx-based \ctlos\ algorithmic outline}
    % : An overview of the SCP Algorithm.
    \label{fig:scp}
    \vspace{-10pt}
\end{figure}

In this section, we will present a brief problem-specific overview of the constraint reformulation, time dilation, control parameterization, linearization, discretization, and convex subproblem steps of the \ctscvx\ method. 
The overall algorithm is visualized in Fig.~\ref{fig:scp}. 
A comprehensive overview written in a more general context of the \ctscvx\ method can be found in \cite{Elango2024ctcs}.
\subsection{Isoperimetric Constraint Reformulation}\label{sec:ctcs}
To ensure continuous-time constraint satisfaction, \eqref{eq:los_constr} is reformulated as
\begin{multline}
    \label{eq:ctcs_reform}
    g_\mathrm{LoS}(t,x(t),u(t)) \leq 0 \ \forall t \in [t_i,t_f] \Longleftrightarrow \cr \int^{t_f}_{t_i}\max\{0, g_\mathrm{LoS}(t,x(t),u(t))\}^2\,dt = 0.
\end{multline}
To satisfy \eqref{eq:los_constr}, it is sufficient to satisfy the RHS of \eqref{eq:ctcs_reform}, which was proposed and proved in Lemma 2 and its proof in \cite{Elango2024ctcs}. A new state is defined as
\begin{align*}
    \dot{y}(t) &= \max\{0, g_\mathrm{LoS}(t,x(t),u(t))\}^2,\\
    y(t_i) &= y(t_f).
\end{align*}
The new state, $y:\mathbb{R}_+ \mapsto \mathbb{R}$, is concatenated to the original state vector forming a new augmented state vector
\begin{align}
    \label{eq:aug_state}
    \tilde{x}(t) &= \begin{bmatrix} x(t)^\top & y(t)^\top \end{bmatrix}^\top,
\end{align}
and augmented dynamics
\begin{align}
    \label{eq:aug_dyn}
    \dot{\tilde{x}}(t) &= \begin{bmatrix} f_{\mathrm{6DOF}}(t,x(t), u(t)) \\ \max\{0, g_\mathrm{LoS}(t,x(t),u(t))\}^2 \end{bmatrix}.
\end{align}

\subsection{Time Dilation}
Time dilation is shown as an effective method to solve free-final-time problems with a uniform time grid in \cite{Szmuk2018free} and with an adaptive time grid in \cite{Kamath2023free}. 
Time dilation reformulates a free-final-time problem as an equivalent fixed-final-time problem by shrinking or expanding, hence ``dilating", the time interval and treating the dilation constant, $s\in\mathbb{R}$, as a control input. 
To allow each discrete node to freely move in the time-domain, we adopt the time dilation with an adaptive time grid. 
A uniform normalized time grid, $\tau \in [0,1]$, is introduced, and the mapping back to the original time interval is $t(\tau) : [0,1] \to [t_i, t_f]$. 
The time-dilation constant is defined as
\begin{align*}
    s(\tau) \triangleq \frac{dt(\tau)}{d\tau}>0.
\end{align*}
We treat the time dilation constant, $s(\tau)$, as a control input and define the augmented control vector as
\begin{align}
    \label{eq:ctrl_aug}
    \tilde{u}(\tau) = \begin{bmatrix} u(t(\tau))^\top & s(\tau) \end{bmatrix}^\top.
\end{align}
We can express \eqref{eq:aug_dyn} in normalized time
\begin{align}
    \label{eq:tau_dy}
    \begin{split}
    \overset{\circ}{\tilde{x}}(\tau) &= \frac{d\tilde{x}(t(\tau))}{d\tau}\\
    &= \frac{d\tilde{x}(t(\tau))}{dt(\tau)}\frac{dt(\tau)}{d\tau}\\
    &= s(\tau) \begin{bmatrix} f_{\mathrm{6DOF}}(t(\tau),x(t(\tau)), u(t(\tau))) \\ \dot{y}(t(\tau)) \end{bmatrix}\\
    &= F(\tilde{x}(\tau), \tilde{u}(\tau)),
    \end{split}
\end{align}
where the derivative with respect to $\tau$ is denoted as $\overset{\circ}{\Box}$. In the remainder of this paper all controls, $u$, are the augmented form defined in \eqref{eq:ctrl_aug}. We omit the tilde operator on the control for notational simplicity.
\subsection{Control Parameterization}
We use a first-order hold (FOH) instead of zero-order hold for control parameterization as the resulting trajectory and control are smoother and easier to track by real-world vehicles \cite{Szmuk2019-bt}. For discrete control signals, $u_{k}$ and $u_{k+1}$, the continuous control signal on the interval $[\tau_k, \tau_{k+1})$ is defined as follows
\begin{equation*}
u(\tau) \triangleq \sigma^{-}_k(\tau)u(k) + \sigma^{+}_{k+1}(\tau)u(k+1), \ \forall \tau \in [\tau_k, \tau_{k+1}),
\end{equation*}
where
\begin{align*}
\sigma_{k}^{-}(\tau) \triangleq \frac{\tau_{k+1}-\tau}{\tau_{k+1}-\tau_k}, \sigma_{k+1}^+(\tau) \triangleq \frac{\tau-\tau_k}{\tau_{k+1}-\tau_k}.
\end{align*}
\subsection{Linearization}\label{sec:lin}
A byproduct of the continuous-time reformulation method proposed in \cite{Elango2024ctcs} is that
the problem is reformulated into a quadratic program. 
To pose our problem in this form, we need to approximate all nonlinear constraints with a first-order Taylor series approximation about a reference trajectory, $(\bar{\tilde{x}}, \bar{u})$.  The Jacobians of the dynamics with respect to $\tilde{x}$ and $u$ are
\begin{align*}
    A(\tau) &\triangleq \nabla_{\tilde{x}}F(\bar{\tilde{x}}(\tau), \bar{u}(\tau)),\\
    B(\tau) &\triangleq \nabla_{u}F(\bar{\tilde{x}}(\tau), \bar{u}(\tau)).
\end{align*}
The linear time-varying (LTV) continuous-time dynamics are
\begin{align}
    \label{eq:lin_dy}
    \Delta \overset{\circ}{\tilde{x}}(\tau) = A(\tau)\Delta\tilde{x}(\tau) + B(\tau)\Delta u(\tau), % + S\Delta \tau,
\end{align}
where $\Delta \Box \triangleq \Box - \bar{\Box}$.

\subsection{Discretization}\label{sec:dis}
To make Problem~\ref{eq:gen_prob} solvable, we discretize it into a finite-dimensional form. 
To maintain dynamic feasibility and avoid constraint violations between discrete points, we use the inverse-free exact discretization method from \cite{Elango2024ctcs}, applied to the dynamics in \eqref{eq:lin_dy}. 
This method is exact, as it relies on the integral form of the differential equation~\eqref{eq:lin_dy}. 

\subsubsection{Discretized Dynamics}
The LTV dynamics with a FOH control parametrization are
\begin{multline}
    \label{eq:dis_dy}
    \Delta \overset{\circ}{\tilde{x}}(\tau) \approx A(\tau)\Delta \tilde{x}(\tau) + B(\tau)\sigma^-_k(\tau)\Delta u(\tau_{k}) + \cr B(\tau)\sigma^+_k(\tau)\Delta u(\tau_{k+1}),
\end{multline}
for $\tau\in[\tau_k, \tau_{k+1})$.
The unique solution to \eqref{eq:dis_dy} is given by \eqref{eq:dis_sol}, as shown in \cite{Malyuta2022csm, Antsaklis2005unq}, as
\begin{multline}
    \label{eq:dis_sol}
    \Delta \tilde{x}(\tau) = \Phi_{k}(\tau, \tau_k) \Delta\tilde{x}(\tau_k) + \int^{\tau_{k+1}}_{\tau_k} \Phi_{k}(\tau, \xi) \cr \{B(\xi)\sigma^-_k(\tau)\Delta u_k + B(\xi)\sigma^+_k(\tau)\Delta u_{k+1}\}\, d\xi.
\end{multline}
The state-transition matrix associated with Equation~\eqref{eq:dis_sol}, denoted by $\Phi_{k}(\tau, \tau_k), \tau\in[\tau_k,\tau_{k+1})$ satisfies the differential equation
\begin{align*}
    \overset{\circ}{\Phi}_{k}(\tau, \tau_k) =  A(\tau)\Phi_{k}(\tau,\tau_k), \quad \Phi_{k}(\tau_k,\tau_k) = \mathbf{I}_{n_{\tilde{x}}}.
\end{align*}
The LTV discretized dynamics are
\begin{subequations}
    \begin{align}
        \Delta\tilde{x}_{k+1} &= \bar{A}_k\Delta\tilde{x}_k + \bar{B}_k^-\Delta u_k + \bar{B}_k^+\Delta u_{k+1},\\
        \bar{A}_k &\triangleq \Phi_{k}(\tau_{k+1}, \tau_k), \label{eq:a_dis}\\
        \bar{B}^{-}_{k} &\triangleq \int^{\tau_{k+1}}_{\tau_k}\Phi_{k}(\tau_{k+1}, \xi)B(\xi)\sigma_{k}^-(\xi)\, d\xi, \label{eq:b_min_dis}\\
        \bar{B}^+_k &\triangleq \int^{\tau_{k+1}}_{\tau_k}\Phi_{k}(\tau_{k+1}, \xi)B(\xi)\sigma_{k}^+(\xi)\, d\xi\label{eq:b_plus_dis}.
    \end{align}
\end{subequations}

\subsection{Convex Subproblem}\label{sec:cvx_sub}

To formulate the convex subproblem, it is necessary to address numerical stability concerns as well as employ a trust region to ensure convergence.

\label{sec:scale}
For numerical stability, we adopt an affine scaling method on the state and control similar to the method presented in \cite{Reynolds2020scaling}. The scaled variables are denoted by the $\hat{\Box}$ operator.

A soft trust region is penalized to ensure the solution, $(x,u)$, does not deviate too far from the reference or previous solution, $(\bar{x},\bar{u})$. 
Specifically, we use a squared 2-norm, $\lVert\Box - \bar{\Box}\rVert_2^2$, on the state and control, which has convergence guarantees (see Section 5 in \cite{Drusvyatskiy2018-pn}). 
\subsubsection{Min-Max State and Control Constraints}
\label{sec:min-max}
To ensure the convergence of Algorithm~\ref{alg:ptr}, the state and control sets, $\mathcal{X}\ \text{and}\ \mathcal{U}$ respectively, are assumed to be compact \cite{Malyuta2022csm}. 
In practice, this is enforced by defining minimum and maximum bounds for each element of the state and control vectors
\begin{subequations}
    \label{eq:min_max_constr}
    \begin{align}
        x^{i}_{\min} &\leq x^{i}_{k} \leq x^{i}_{\max}, &&\forall i \in \{0,n_{x}\}, \forall k \in\{0,N\},\\
        u^{i}_{\min} &\leq u^{i}_{k} \leq u^{i}_{\max}, &&\forall i \in \{0,n_{u}\}, \forall k \in\{0,N\}.
    \end{align}
\end{subequations}
The boundary condition on the augmented state $y$ is relaxed to
\begin{align*}
    y_k - y_{k-1} &\leq \varepsilon_{\mathrm{LICQ}} && \forall k \in \{1,N\},
\end{align*}
where $\varepsilon_{\mathrm{LICQ}} \geq 0$ is sufficiently small ($\sim 10^{-4}$) such that all feasible solutions do not violate the linear independence constraint qualification (LICQ) outlined in Section 3.1 of \cite{Elango2024ctcs}. 
This ensures that the gradients of active constraints are linearly independent.
The constraints in \eqref{eq:min_max_constr} are enforced using the same continuous-time constraint reformulation method outlined in \ref{sec:ctcs}. 
\subsubsection{Convex Subproblem Form}

The subproblem is formed by applying the  aforementioned constraint reformulation, discretization, and linearization steps to the original NLP \eqref{eq:gen_prob}.
\label{sec:cvx_sub_form}

\begin{mybox}{\textbf{CT-LoS Convex Subproblem}}
\small
\begin{subequations}
\label{eq:ctlos_sub}
\begin{alignat}{2}
    &\!\minimize{\hat{\tilde{x}},\hat{u}} &\quad& \lambda_{\mathrm{obj}}L_{f}(\tilde{x}_{N},u_{N}) +\label{eq:1} \\
         & & & \sum^N_{k=0}  \lambda_{\mathrm{tr}}\left|\left|\begin{bmatrix} \hat{\tilde{x}}_{k} \\ \hat{u}_{k}\end{bmatrix} - \begin{bmatrix} \bar{\hat{\tilde{x}}}_k \\ \bar{\hat{u}}_{k} \end{bmatrix}\right|\right|_2^2 + \lambda_{\mathrm{vc}}\lVert\nu_k\rVert_1\nonumber\\
    &\st & &\Delta\tilde{x}_{k} = \bar{A}_{k-1}\Delta\tilde{x}_{k-1} + \bar{B}_{k-1}^-\Delta u_{k-1}\nonumber \\ 
        & & &\quad + \bar{B}_{k-1}^+\Delta u_{k} + \nu_{k-1}, \\
    & & & y_k - y_{k-1} \leq \varepsilon_{\mathrm{LICQ}}, \label{eq:licq}\\
    & & &\tilde{x}_0 = \tilde{x}_i,\tilde{x}_N = \tilde{x}_f,
\end{alignat}
\end{subequations}
where $k \in \{1,\dots, N\}$.
\end{mybox}
\normalsize
\noindent The objective weight $\lambda_{\Box}$ corresponds to the $\Box$ term.

\subsection{Prox-Linear Method}
The prox-linear method is outlined in Alg.~\ref{alg:ptr}.

\begin{figure}[h]
\begin{myboxyellow}\vspace{-4mm}
\begin{minipage}{\linewidth}
\begin{algorithm}[H]
\caption{Prox-Linear Method}
\small
\label{alg:ptr}
\begin{algorithmic}[0]
\Require $\epsilon_{\mathrm{tr}}$, $\epsilon_{\mathrm{vc}}, \epsilon_{\mathrm{vb}}$, $\bar{x}, \bar{u}, k_{\max}$
\State \textbf{Initialize} $x \ne \bar{x}, u \ne \bar{u}$, $\xi_{\mathrm{vc}} \ne 0$, $\xi_{\mathrm{vb}} \ne 0$, $k=1$
\While{$k < k_{\max}\textbf{ and }$ $G(x,u,\bar{x}, \bar{u}, \xi_{\mathrm{vc}},  \xi_{\mathrm{vb}})$}
\State $x,u,\xi_{\mathrm{vc}},\xi_{\mathrm{vb}}\gets$ Solve Convex Subproblem, either \eqref{eq:ctlos_sub} or \eqref{eq:dtlos_sub}, with linearization about $\bar{x},\bar{u}$
\State $\bar{x},\bar{u} \gets x,u$ 
\State $k \gets k+1$
\EndWhile\\
\Return $x,u$
\end{algorithmic}
\end{algorithm}
\end{minipage}\vspace{2mm}
where,
\begin{multline*}
    G(x,u,\bar{x}, \bar{u}, \xi_{\mathrm{vc}},  \xi_{\mathrm{vb}}) = \left\lVert \begin{bmatrix}
    x \\ u\end{bmatrix} - \begin{bmatrix} \bar{x} \\ \bar{u}\end{bmatrix}\right\rVert_2^2 > \epsilon_{\mathrm{tr}}\cr\textbf{ or }\lVert\xi_{\mathrm{vc}}\rVert_1 > \epsilon_{\mathrm{vc}} \textbf{ or } \lVert \xi_{\mathrm{vb}}\rVert_{1} > \epsilon_{\mathrm{vb}},
\end{multline*}
is a boolean operator to determine convergence.
\end{myboxyellow}
\end{figure}
\normalsize
\noindent The virtual buffer, $\xi_{\mathrm{vb}}$, models linearized nonconvex path constraint violations at discrete nodes, while the virtual control, $\xi_{\mathrm{vc}}$, models continuous-time dynamic violations. In the \ctlos\ formulation, linearized path constraints are appended to the dynamics, eliminating the need for a virtual buffer as these violations are now captured in the virtual control.

The convergence criterion for the trust region, and virtual control and buffer are $\epsilon_{\mathrm{tr}}, \epsilon_{\mathrm{vc}},$ and $\epsilon_{\mathrm{vb}}$. The interested reader may refer to Section 4.2 in \cite{Elango2024ctcs} for additional details regarding the prox-linear method.

\section{Experiments}\label{sec:experiments}

We evaluate our approach on two representative LoS guidance problems and address the following questions:
\begin{itemize}
    \item[\textbf{Q1.}] How well does \ctlos\ satisfy the LoS constraint throughout the trajectory compared to the baseline?
    \item[\textbf{Q2.}] What tradeoffs are made to achieve better LoS violation performance?
    \item[\textbf{Q3.}] How does \ctlos\ scale as the problem size increase?
\end{itemize}
To address the above guiding questions, we will use the following metrics to judge the performance of \ctlos\ and \dtlos\ across both experiments: 1. LoS constraint violation over the trajectory, $\mathrm{LoS}_{\mathrm{vio}}$, which is defined as
    \begin{align*}
        \mathrm{LoS}_{\mathrm{vio}} &\triangleq\frac{\sum_{i=0}^{N_{\mathrm{prop}}}\max\{0,g_{\mathrm{LoS}}(x^i,u^i)\}}{N_{\mathrm{prop}}},
    \end{align*}
where $N_{\mathrm{prop}}$ is the number of propagated nodes $(\sim 1000)$, 2. the runtime and number of iterations till convergence of Alg.~\ref{alg:ptr}, and 3. the original objective cost (total fuel consumption in the cinematography scenario, time of flight in the relative navigation scenario).

We systematically vary the discretization grid size to adjust the problem size and observe the impact on metrics, while also sweeping through objective weights for each grid size to avoid over-tuning and ensure the problems are not cherry-picked.

\subsection{Baseline Comparison}
\label{sec:dt_comp}
We present a baseline comparison to the \ctlos\ formulation, which enforces LoS constraint violation at discrete nodes rather than the integral of LoS constraint violation and refer to it as the \dtlos\ formulation. Except for the LoS constraint, the \dtlos\ formulation follows the same formulation and solution steps as the SCvx methods presented in \cite{Malyuta2022csm}.
\subsubsection{Discrete-Time Constraints}
All nonconvex constraints are linearized 
\begin{align}
    \label{eq:lin_constr}
    \begin{split}
         L_{\Box}(x_k, u_k) &= \Box(\bar{x}_k, \bar{u}_k) + \\ \nabla_x \Box(\bar{x}_k, &\bar{u}_k)(x_k-\bar{x_k}) + 
        \nabla_{u} \Box(\bar{x}_k, \bar{u}_k)(u_k-\bar{u}_k),
    \end{split}
\end{align}
where $\Box = g$ for inequality path constraints and $\Box = h$ for equality path constraints. 
All linearized constraints are enforced using a slack variable and all convex constraints are directly enforced at each node.
\begin{mybox}{\textbf{DT-LoS (Baseline) Convex Subproblem}}
\small
\begin{subequations}
\label{eq:dtlos_sub}
\begin{alignat}{2}
    &\!\minimize{\hat{x},\hat{u}} &\quad& \lambda_{\mathrm{obj}}L_{f}(x_{N},u_{N}) 
 + \lambda_{\mathrm{vc}}\lVert\nu_k\rVert_1 + \nonumber\\
 & & & + \lambda_{\mathrm{vb}}\max\{0,\nu_{g_{\mathrm{LoS}}}\} \label{eq:1} \\
         & & & +  \sum^N_{k=0}  \lambda_{\mathrm{tr}}\left|\left|\begin{bmatrix} \hat{x}_k \\ \hat{u}_{k}\end{bmatrix} - \begin{bmatrix} \bar{\hat{x}}_k \\ \bar{\hat{u}}_{k} \end{bmatrix}\right|\right|_2^2 \nonumber\\
    &\st & &\Delta x_{k} = \bar{A}_{k-1}\Delta x_{k-1} + \bar{B}_{k-1}^-\Delta u_{k-1}\nonumber \\ 
        & & &\quad + \bar{B}_{k-1}^+\Delta u_{k} + \nu_{k-1}, \\
    & & &L_{g_{\mathrm{LoS}}}^{\ell}(x_{k}, u_{k}) = \nu_{g_{\mathrm{LoS}}}^{\ell}, \\
    & & &x^{i}_{\min} \leq x^{i}_{k} \leq x^{i}_{\max}, \\
    & & &u^{j}_{\min} \leq u^{j}_{k-1} \leq u^{j}_{\max},\\
    & & &x_0 = x_i, \ x_N = x_f,
\end{alignat}
\end{subequations}
where appropriate, $\forall i \in \{0,\dots, n_{x}\}, \forall j \in \{0,\dots, n_{u}\}, \forall k \in \{1,\dots, N\}, \forall \ell \in \{1, \dots, n_{kp}\}$.
\end{mybox}
\normalsize

\subsection{Cinematography Scenario}

\begin{figure}[t] 
    \centering
    \begin{subfigure}[t]{0.9\linewidth}
        \centering
        \includegraphics[width=0.98\linewidth]{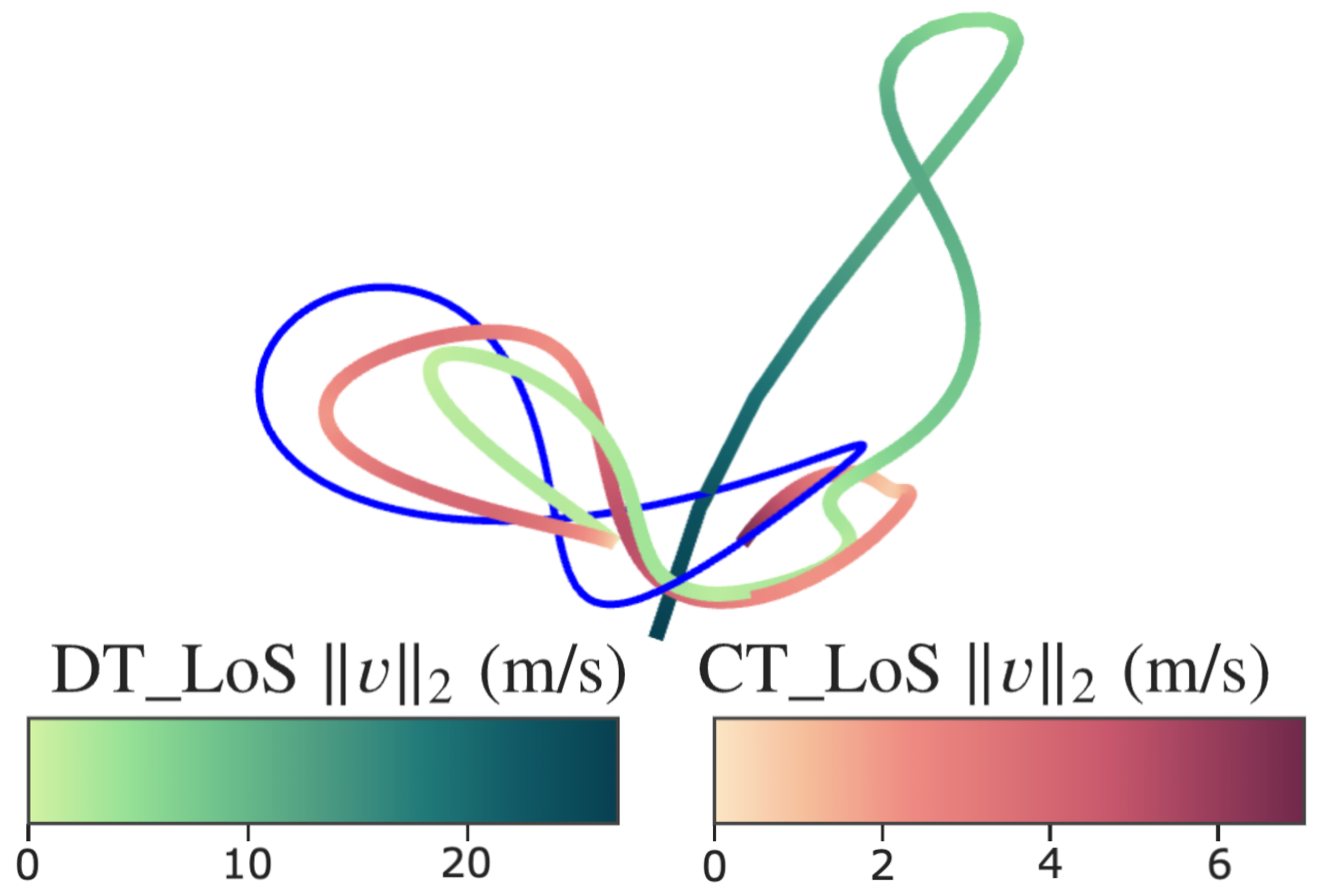}
        \caption{\footnotesize Cinematography Scenario.}
        \label{fig:cine_sce}
    \end{subfigure}
    % \hspace{0.02\linewidth}
    \begin{subfigure}[t]{0.9\linewidth}
        \centering
        \includegraphics[width=0.98\linewidth]{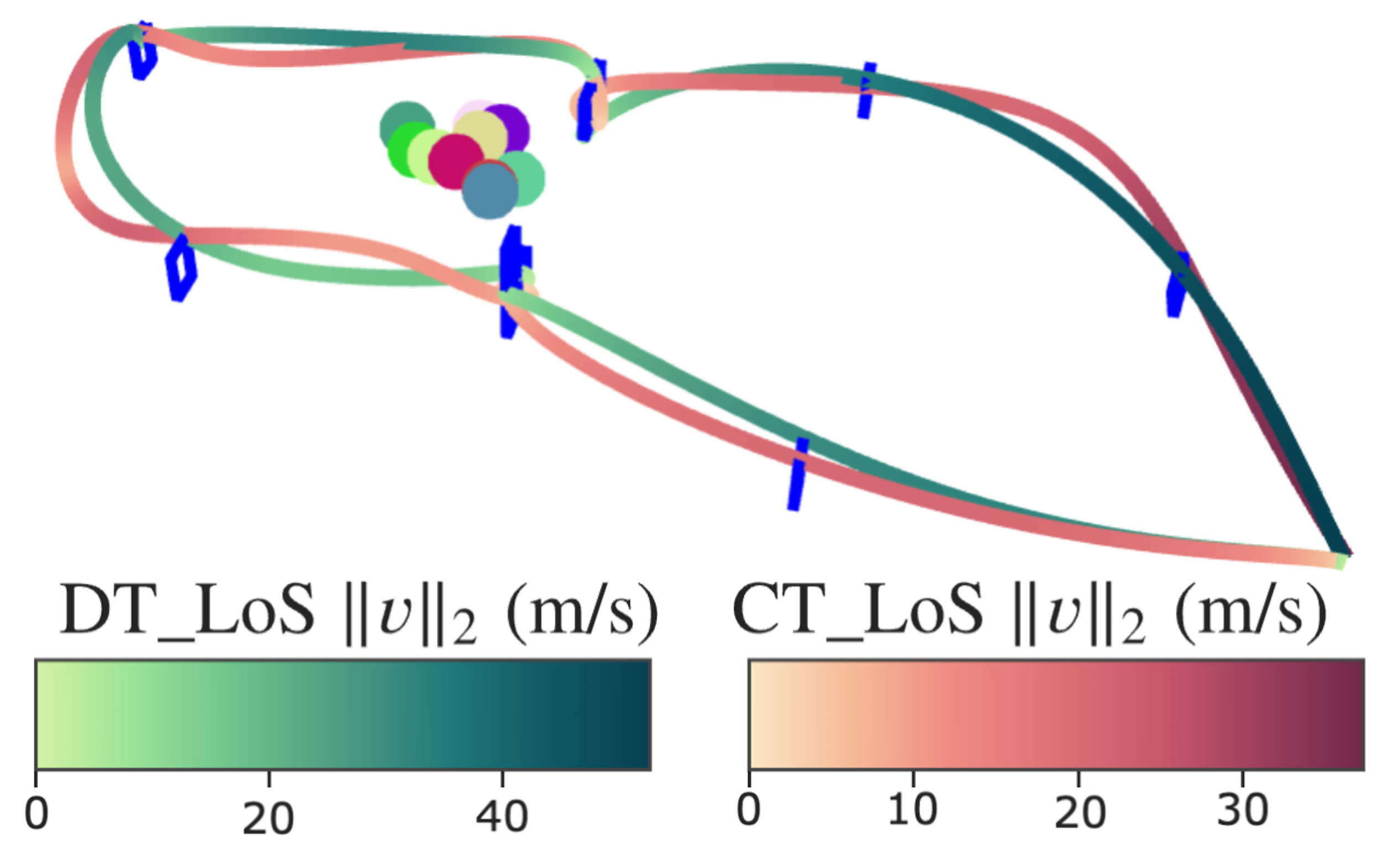}
        \caption{\footnotesize Relative Navigation Scenario}
        \label{fig:dr_sce}
    \end{subfigure}
    \caption{Cinematography and Relative Navigation Scenarios: The drone’s path is shown with a color gradient indicating the $L_2$-norm of velocity for both the baseline and proposed methods. The keypoints are denoted by dots. In the cinematography scenario, the keypoints’ trajectory is shown in solid blue.}
    \label{fig:combined_scenarios}
    \vspace{-10pt}
\end{figure}

The cinematography scenario, seen in Fig.~\ref{fig:cine_sce}, draws inspiration from the LoS guidance challenge in drone-enabled cinematography. In this scenario, a single dynamic keypoint with a known trajectory, representing the subject to be filmed, must remain within LoS of the drone. Additionally, the drone must be within a minimum and maximum range of the subject. The subject’s trajectory follows a simple sinusoidal function. A non-symmetric $\ell_{\infty}$-norm cone is used to model the viewcone. Since the terminal-state of the subject implicitly determines the time-of-flight (ToF) and terminal-state constraint, the ToF is arbitrarily fixed, and there is no explicit terminal-state constraint. The objective of this scenario is to minimize fuel cost, $\int^{t_f}_{t_0}||u(t)||_2\,dt.$

\subsection{Relative Navigation Scenario}
\begin{figure}[t]
    \centering
    \includegraphics[width=\linewidth]{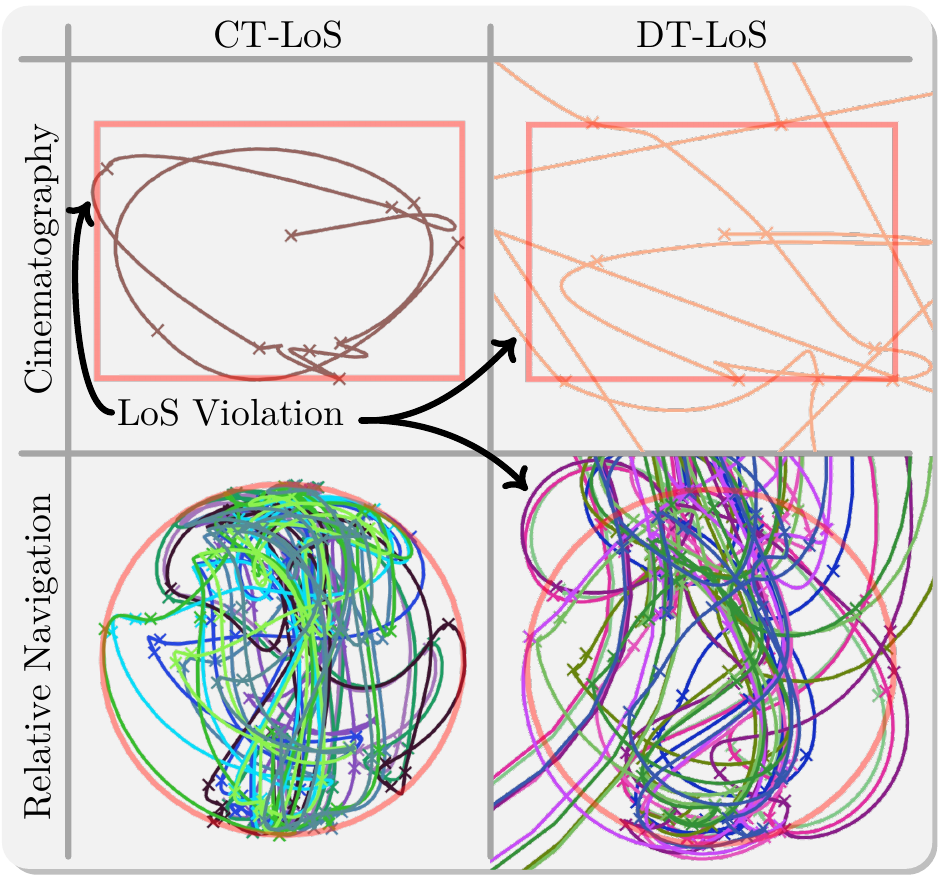}
    \caption{Qualitative Results: The trajectories of keypoints in the sensor frame for both \ctlos\ and \dtlos\ methods are displayed. 
    $\times$'s represent discrete nodes, while solid lines represent nonlinear keypoint paths. Different colors are used to distinguish multiple keypoint paths. The solid red lines represent the sensor field-of-view.}
    \label{fig:qual_res}
    \vspace{-10pt}
\end{figure}

The relative navigation experiment, seen in Fig.~\ref{fig:dr_sce}, is inspired by the LoS guidance challenge encountered in vision-aided drone racing \cite{Song2023-mt}. In this scenario, ten static keypoints serve as landmarks for visual-inertial odometry and must remain within the LoS as the drone navigates through ten square gates in a predefined sequence. The goal is to complete this sequence in minimal time. The scenario is visualized in Fig.~\ref{fig:qual_res}. A symmetric $\ell_2$-norm cone is used to model the viewcone. 

\subsection{Implementation Details}\label{sec:imp_det}
The code used to generate the results in this paper is written in Python and can be found here, \href{https://haynec.github.io/papers/los}{https://haynec.github.io/papers/los}. 
During the discretization step, equations \eqref{eq:a_dis},  \eqref{eq:b_min_dis}, and \eqref{eq:b_plus_dis} are evaluated using the 4th-order Runga-Kutta method from the \texttt{SciPy} package, \cite{2020SciPy-NMeth}.
The convex subproblems for both methods are solved using the package \texttt{cvxpy}, \cite{agrawal2018rewriting, diamond2016cvxpy}, using the convex solver, Clarabel \cite{Goulart2024-bz}. 
The Jacobians are calculated using automatic differentiation in \texttt{JAX}, \cite{jax2018github}.
All experiments are run on a workstation with an AMD 7950X3D CPU and 128 GB of memory.

% \subsection{Results}\label{sec:res}

\begin{figure*}[h]
    \centering
    \begin{subfigure}[t]{0.45\linewidth}
        \centering
        \includegraphics[width=\linewidth]{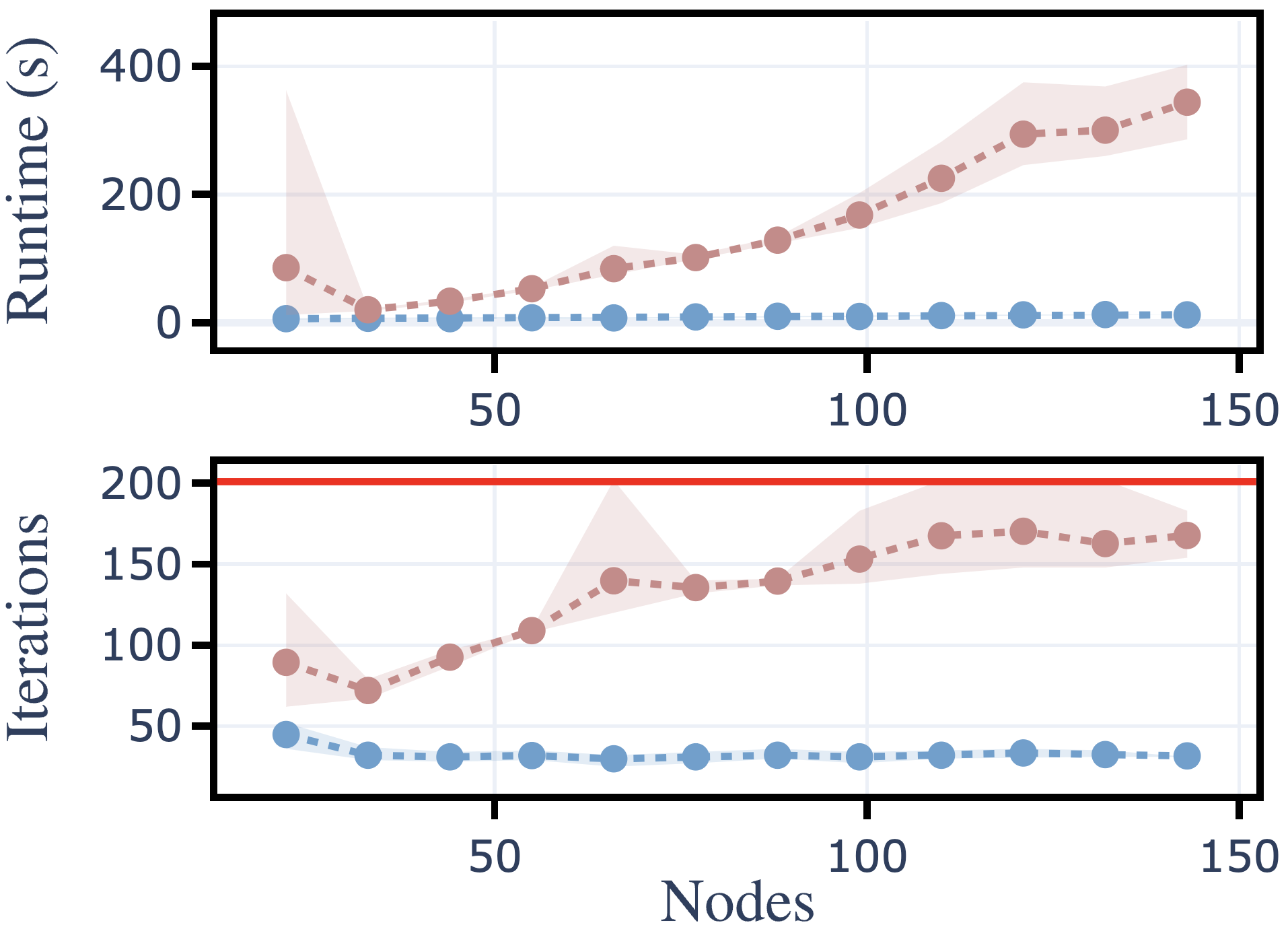}
        \caption{\footnotesize Relative Navigation Runtime}
        \label{fig:dr_run}
    \end{subfigure}
    \hfill
    \begin{subfigure}[t]{0.45\linewidth}
        \centering
        \includegraphics[width=\linewidth]{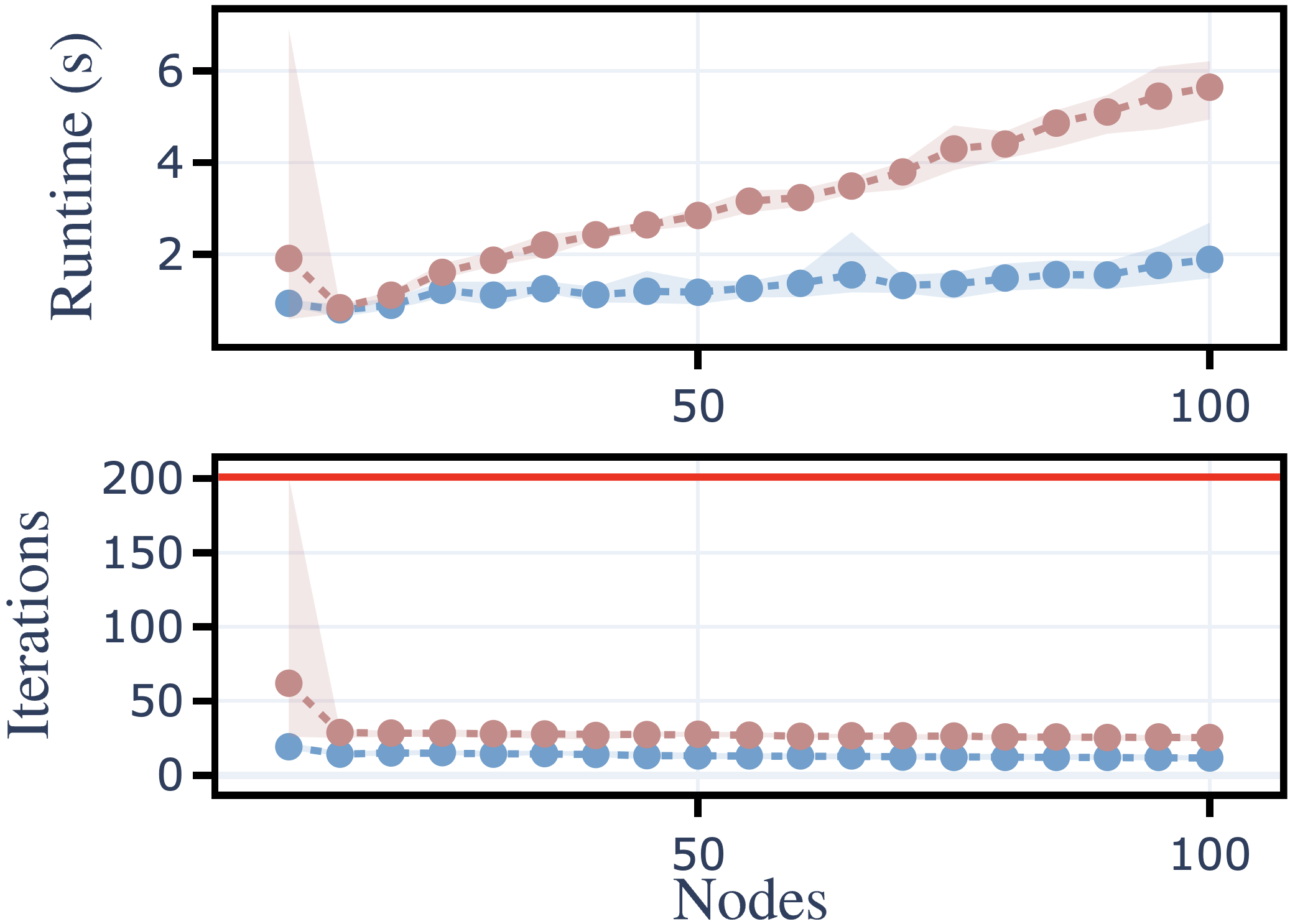}
        \caption{\footnotesize Cinematography Runtime}
        \label{fig:cine_run}
    \end{subfigure}
    \hfill
    \begin{subfigure}[t]{0.45\linewidth}
        \centering
        \includegraphics[width=\linewidth]{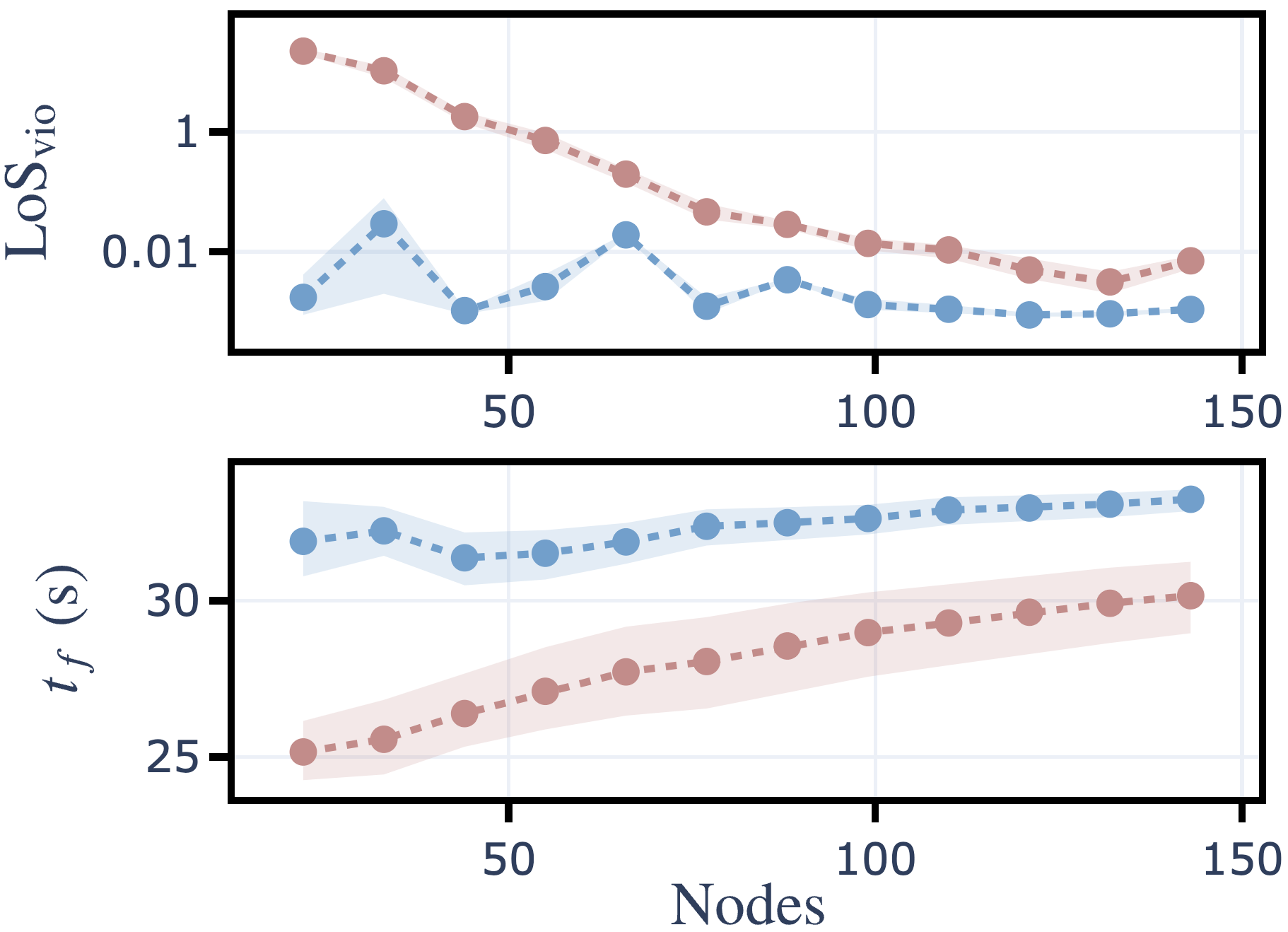}
        \caption{\footnotesize Relative Navigation LoS Violation}
        \label{fig:dr_los}
    \end{subfigure}
    \hfill
    \begin{subfigure}[t]{0.45\linewidth}
        \centering
        \includegraphics[width=\linewidth]{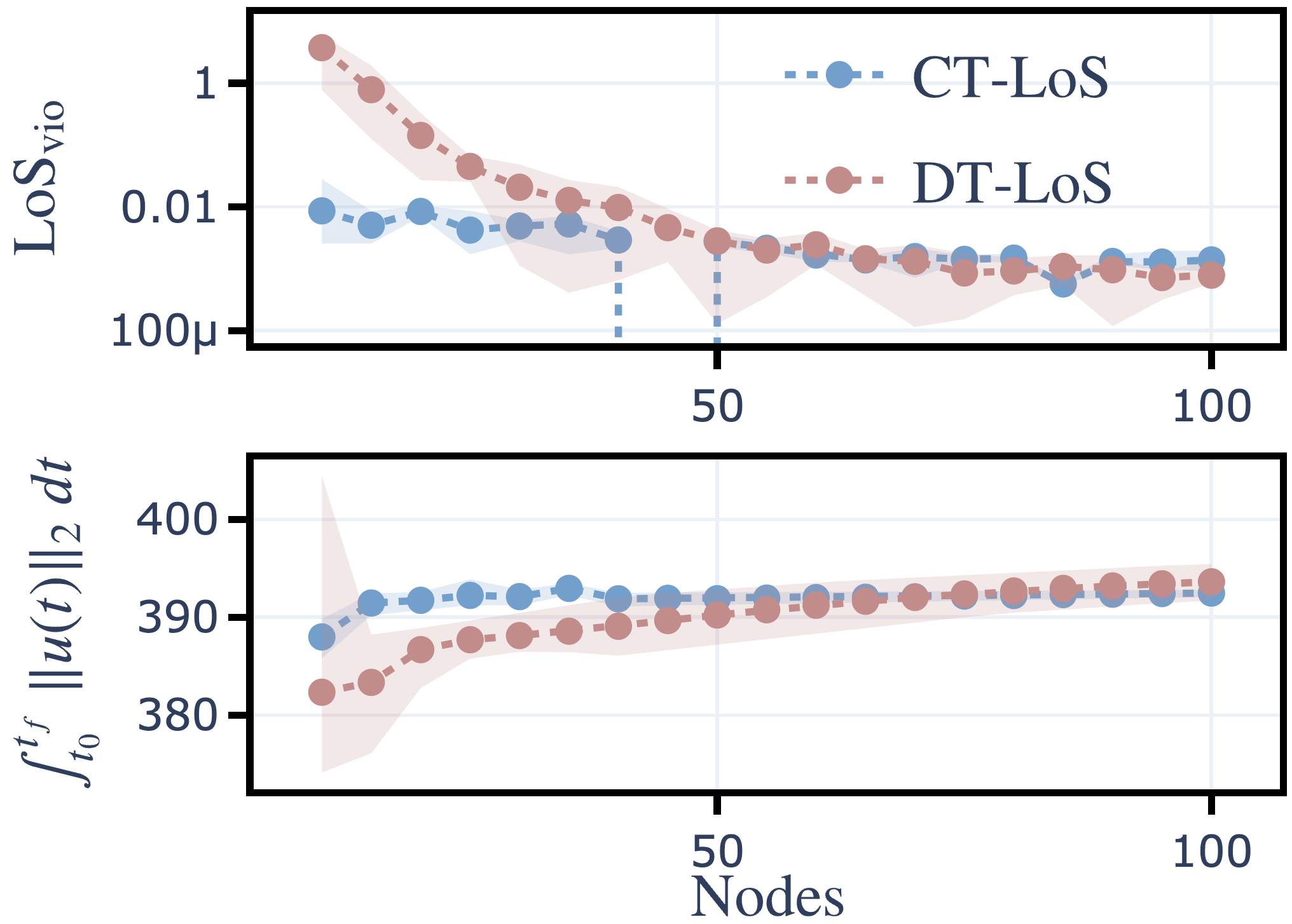}
        \caption{\footnotesize Cinematography LoS Violation}
        \label{fig:cine_los}
    \end{subfigure}
    \caption{\small Quantitative Results: (a) and (b) show the total runtime and number of iterations until convergence of Alg.~\ref{alg:ptr} respectively. (c) and (d) show the LoS constraint violation and the time-of-flight or total fuel cost respectively. The shaded regions represent the minimum and maximum and the dots represent the average values of runtime and iterations for the upper and lower plots respectively. For all plots, lower is better. The maximum number of iterations is denoted by the red horizontal line in the iteration plots.}
    \label{fig:combined_grid}
    \vspace{-10pt}
\end{figure*}

\section{Results \& Discussion}\label{sec:discuss}
To address the question of LoS violation (\textbf{Q1}), the relative navigation scenario with a grid size of 22 nodes, \ctlos\ achieves an average $\mathrm{LoS}_{\mathrm{vio}}$ of $1.73\times10^{-3}$, significantly lower than \dtlos\ at $22.35$, as seen in Fig.~\ref{fig:dr_los}. This indicates that \ctlos\ maintains the LoS constraint more effectively throughout the entire trajectory. Similarly, in the cinematography scenario with a grid size of $10$ nodes, \ctlos\ achieves an average $\mathrm{LoS}_{\mathrm{vio}}$ of $8.63\times10^{-3}$, compared to \dtlos\ at $3.76$, as seen in Fig.\ref{fig:cine_los}. These results demonstrate that \ctlos\ consistently either outperforms or has equivalent performance to \dtlos\ in maintaining LoS constraints for all node counts.

Regarding the tradeoffs made to achieve better LoS violation performance (\textbf{Q2}), the primary tradeoff of the proposed approach is the potential sacrifice in objective performance. While \ctlos\ penalizes the integral of $\mathrm{LoS}_{\mathrm{vio}}$ along the trajectory, \dtlos\ only penalizes $\mathrm{LoS}_{\mathrm{vio}}$ at discrete points. This can be seen in Fig.~\ref{fig:cine_los} and Fig.~\ref{fig:dr_los}, in which the \ctlos\ curve is below the \dtlos\ except for large node counts in the cinematography scenario. This is expected as \dtlos\ solves a less constrained approximation of the true NLP we are trying to solve. \ctlos\ directly solves the true NLP by penalizing the integral of $\mathrm{LoS}_{\mathrm{vio}}$ along the trajectory resulting in a more constrained problem.

Concerning the scalability of \ctlos\ as the problem size increases (\textbf{Q3}), runtime for both formulations scales linearly with grid size, but \dtlos\ has a steeper slope, especially in the relative navigation scenario. 
For example, increasing the grid size to $132$ nodes improves \dtlos\ $\mathrm{LoS}_{\mathrm{vio}}$ to $3.14\times10^{-3}$, but at the cost of a prohibitively longer runtime of $3.00\times10^2$ seconds, compared to \ctlos\ at $6.05$ seconds, as seen in Fig.~\ref{fig:dr_los} and Fig.~\ref{fig:dr_run}. In the cinematography scenario, increasing the grid size to $45$ nodes reduces \dtlos\ $\mathrm{LoS}_{\mathrm{vio}}$ to $4.62\times10^{-3}$, but increases runtime to $2.65$ seconds, compared to \ctlos\ at $1.19$ seconds, as seen in Fig.~\ref{fig:cine_los} and Fig.~\ref{fig:cine_run}. This is due to the increased number of constraints added to the \dtlos\ convex subproblem with each additional node. In contrast, \ctlos\ adds only a single LICQ constraint and an additional row to the dynamics constraint per additional node, making it more computationally efficient. Furthermore, \dtlos\ requires at least twice the number of iterations to converge compared to \ctlos, highlighting the scalability advantage of \ctlos\ as the problem size increases.

The \ctlos\ algorithm demonstrates its efficacy across two challenging scenarios inspired by real-world problems. 
Across both scenarios, our algorithm consistently shows either lower or equivalent $\mathrm{LoS}_{\mathrm{vio}}$ than \dtlos. As the grid size increases, \dtlos\ performance converges to that of \ctlos\ due to more frequent enforcement of nodal LoS constraints. The \ctlos\ formulation remains largely unaffected by grid size due to its integral constraint violation penalty. Notably, \ctlos\ does sacrifice objective performance, which can be seen in Fig.~\ref{fig:cine_los} and Fig.~\ref{fig:dr_los}, except for larger grid sizes in the cinematography scenario. 

\section{Conclusion}\label{sec:conclusion}
In this work, we addressed three central challenges within the LoS guidance problem. We proposed a sensor-footprint-agnostic LoS constraint with a convex norm cone and nonconvex transformation components. Additionally, we developed a computationally tractable formulation of a continuous-time LoS guidance method, \ctlos, which leverages the \ctscvx\ algorithm from \cite{Elango2024ctcs}. We then demonstrated the proposed method’s efficacy in several representative and challenging scenarios inspired by real-world applications. Finally, we compared it against a baseline formulation.
This work is a step forward to building reliable real-time planning algorithms suitable for safety-critical applications.

\bibliographystyle{ieeetr}
\bibliography{references}

\end{document}